\documentclass[11pt]{article}
\usepackage{epsfig,makeidx,amsmath,amsfonts,latexsym,subfigure}
\newtheorem{theorem}{Theorem}[section]
\newtheorem{prop}{Proposition}[section]
\newtheorem{lemma}{Lemma}[section]

\newcommand{\E}{{\mathbb E}}
\newcommand{\RR}{{\mathbb R}}
\newcommand{\Z}{{\mathbb Z}}
\newcommand {\PP}{\mathbb P}

\begin{document}

\title{Stationary random graphs with prescribed iid degrees on a spatial Poisson process}

\author{Maria Deijfen \thanks{Department of Mathematics, Stockholm University, 106 91 Stockholm. Email: mia@math.su.se.}}

\date{October 2008}

\maketitle

\thispagestyle{empty}

\begin{abstract}

\noindent Let $[\mathcal{P}]$ be the points of a Poisson process on $\RR^d$ and $F$ a probability distribution with support on the non-negative integers. Models are formulated for generating translation invariant random graphs with vertex set $[\mathcal{P}]$ and iid vertex degrees with distribution $F$, and the length of the edges is analyzed. The main result is that finite mean for the total edge length per vertex is possible if and only if $F$ has finite moment of order $(d+1)/d$.

\vspace{0.5cm}


\vspace{0.5cm}

\end{abstract}

\section{Introduction}

Consider a Poisson process $\mathcal{P}$ on $\mathbb{R}^d$ with intensity 1 and write $[\mathcal{P}]$ for the point set of the process. Furthermore, take a probability distribution $F$ with support on $\mathbb{N}$. How should one do to obtain a translation invariant random graph with vertex set $[\mathcal{P}]$ and degree distribution $F$? Which properties does the resulting configuration have? These are the questions that will be considered, and partly answered, in this paper.

The problem of generating random graphs with prescribed degree distribution has been extensively studied the last few years. It is primarily motivated by the increasing use of random graphs as models for complex networks; see e.g.\ Durrett (2006). Such networks typically exhibit non-Poissonian distributions for the degrees, predominantly various types of power-laws, and it is hence interesting to develop graph models that can produce this type of degree distributions (indeed, the standard Erd\"{o}s-Renyi graph gives Poisson degrees in the limit). One example of a graph model for obtaining a prescribed degree distribution is the so called configuration model -- see Molloy and Reed (1995,1998) -- where, given $n$ vertices, each vertex is independently assigned a random number of stubs according to the prescribed distribution and the stubs are then paired randomly to create edges. See also Bollob\'{a}s et al.\ (2001), Chung and Lu (2001:1,2) and Bollob\'{a}s et al.\ (2006) for other work in the area.

Most existing models for obtaining graphs with given degree distribution do not take spatial aspects into account, that is, there is no metric defined on the vertex set. In Deijfen and Meester (2006) however, a model is introduced where the vertex set is taken to be $\Z$. The vertices are equipped with stubs according to a prescribed distribution and these stubs are then paired according to a certain rule. This model, which turns out to be closely related to random walks, is shown to lead to well-defined configurations, but the edge length has infinite mean. In Deijfen and Jonasson (2006), a different model (on $\Z$) is formulated and shown to give finite mean for the edges. This work is generalized in Jonasson (2007) to more general (deterministic) vertex structures, in particular it is shown there that there exists a translation invariant graph on $\Z^d$ with finite total edge length per vertex if and only if the degree distribution has finite mean of order $(d+1)/d$.

Here we will take the vertex set to be the random point set $[\mathcal{P}]$ of the Poisson process $\mathcal{P}$. Each vertex is independently assigned a number of stubs according to the desired degree distribution $F$ and we then ask for a way of connecting the stubs to form edges. The restrictions are that the procedure should be translation invariant and the resulting graph is not allowed to contain self-loops and multiple edges, that is, an edge cannot run from a vertex to itself and each pair of vertices can be connected by at most one edge. A few different ways of doing this will be discussed, and we will show that the condition for the possibility of finite mean for the total edge length of a vertex is the same as on $\Z^d$.

As described above, the model can (at least for $d=1,2$) be seen as a way of introducing geography in a complex network model. However, it can also be viewed as a generalization of the problem treated in Holroyd et al.\ (2008), where the matching distance in various types of matchings of the points of a Poisson process is analyzed. Indeed, this would correspond to taking $F\equiv 1$ in our setup, that is, to equip each vertex with one single stub. We will make use of results from Holroyd et al.\ at several occasions.

To formally describe the quantities we will work with, let $(\mathcal{P},\xi)$ be a marked Poisson process, where the ground process $\mathcal{P}$ is a homogeneous Poisson process on $\RR^d$ with rate 1 and the marks are iid with distribution $F$; see Daley and Vere-Jones (2002, Section 6.4). We will often refer to the points in $[\mathcal{P}]$ as vertices, and we think of the marks as representing an assignment of stubs to the vertices. Let $\mathcal{M}$ be a translation invariant scheme for pairing these stubs so that no two stubs at the same vertex are paired and at most one pairing is made between stubs at two given vertices. Furthermore, let $(\mathcal{P}^*,\xi^*,\mathcal{M}^*)$ be the Palm version of $(\mathcal{P},\xi,\mathcal{M})$ with respect to $\mathcal{P}$ and write $\PP^*$ and $\E^*$ for the associated probability law and expectation operator respectively, that is, $\PP^*$ describes the conditional law of $(\mathcal{P},\xi,\mathcal{M})$ given that there is a point at the origin, with the mark process and the pairing scheme taken as stationary background; see Kallenberg (1997, Chapter 11) for details. Note that, since the Palm version of a homogeneous Poisson process has the same distribution as the process itself with an added point at the origin, we have that $[\mathcal{P}^*]\stackrel{d}{=}[\mathcal{P}]\cup \{0\}$. Let $T$ denote the total length in $\mathcal{M}^*$ of all edges at the origin vertex. We will show the following result for $T$.

\begin{theorem}\label{th:main}
There exists a translation invariant pairing scheme with $\E^*[T]<\infty$ if and only if $F$ has finite moment of order $(d+1)/d$.
\end{theorem}

Let $D$ denote the degree of the origin. The only-if direction of the theorem is fairly easy and follows from investigating the radius of the smallest ball around the origin of $\RR^d$ that contains $D$ points in $\mathcal{P}^*$. Indeed, since multiple pairing to the same vertex is not allowed, $T$ is bounded from below by the sum of the distances to the vertices in this ball, which, as it turns out, cannot have finite mean if $\E[D^{(d+1)/d}]$ is infinite. As for the if-direction, we will describe a concrete pairing scheme with the property that $\E^*[T]<\infty$. Here, results from Jonasson (2007) and Holroyd et al. (2008) will be useful.

Before proceeding, we remark that the distribution of $T$ can be described in an alternative way: For $r\in[0,\infty]$, let $H(r)$ denote the expected number of points in $[\mathcal{P}]\cap[0,1]^d$ whose total edge length does not exceed $r$. Since $\mathcal{P}$ has intensity 1 it is not hard to see that $H$ is a distribution function, and, in fact, the distribution of $T$ is given by $H$; see Holroyd et al. (2008, Section 2) for a more careful description of this connection.

The rest of the paper is organized so that the only-if direction of Theorem \ref{th:main} is proved in Section 2, a possible pairing schemes is described in Section 3 and a refinement of this scheme that proves the if-part of Theorem \ref{th:main} is formulated in Section 4. Section 5 contains an outline of possible further work.

\section{Only-if}

Recall that $D$ denotes the degree of the origin vertex in $\mathcal{P}^*$. We will show that $\E^*[T]<\infty$ is impossible if $\E[D^{(d+1)/d}]=\infty$. To this end, let $B(r)$ be the $d$-dimensional ball with radius $r$ centered at the origin and write $V^*(r)$ for the number of $\mathcal{P}^*$-points in $B(r)$, excluding the origin, that is,
$$
V^*(r)=\textrm{card}[\mathcal{P}^*\cap B(r)]-1,
$$
where $\textrm{card}[\cdot]$ denotes set cardinality. Furthermore, for $n\geq 1$, define $R_n$ to be the radius of the smallest ball centered at the origin that contains at least $n$ points of $\mathcal{P}^*$ not counting the point at the origin, that is,
$$
R_n=\inf\{r: V^*(r)\geq n\}.
$$
The idea is that, with a probability that can be bounded from below, at least half of the $D$ stubs at the origin will be connected to stubs at vertices whose distance to the origin is at least $3^{-1/d}R_D$. On this event we have that $T\geq D R_D/(2\cdot 3^{-1/d})$ and, as we will see, the mean of $R_n$ behaves like $n^{1/d}$ for large $n$. Combining this will give that $T$ must have infinite mean if $D^{1+1/d}$ has so.

Consider a Poisson process $\mathcal{P}_\lambda$ on $\RR^d$ with intensity $\lambda$ and -- in analogy with the above definitions -- let $V(r)=\textrm{card}[\mathcal{P}_\lambda\cap B(r)]$ and $R_n=\inf\{r:V(r)\geq n\}$. The mean of $R_n$ is characterized in the following lemma.

\begin{lemma}\label{le:vvRn}
For $n\geq 1$, we have that

\begin{equation}\label{eq:vvRn}
\E[R_n]=\frac{C}{\lambda^{1/d}}\cdot\sum_{k=0}^{n}\frac{\Gamma(k+\frac{1}{d})}{\Gamma(k+1)},
\end{equation}

\noindent where $C=C(d)$ is a constant. In particular, $\E[R_n]\geq C'\left(\frac{n}{\lambda}\right)^{1/d}$ for all $n\geq 1$ and some constant $C'=C'(d)$.
\end{lemma}

\noindent \textbf{Proof.} Since the volume of $B(r)$ is $cr^d$, where $c=c(d)$ is a constant, we have that $V(r)$ is Poisson distributed with parameter $\lambda cr^d$. Together with the definition of $R_n$, this yields that
$$
\PP(R_n\geq r)=\PP(V(r)\leq n)=\sum_{k=0}^n\frac{(\lambda cr^d)^k}{k!}e^{\lambda cr^d}.
$$
Hence
$$
\E[R_n]=\int_0^\infty \PP(R_n\geq r)dr=\sum_{k=0}^n\frac{1}{k!}\int_0^\infty(\lambda cr^d)^ke^{\lambda cr^d}dr,
$$
and, by variable substitution,
$$
\int_0^\infty(\lambda cr^d)^ke^{\lambda cr^d}dr=\frac{\Gamma(k+\frac{1}{d})}{d(\lambda c)^{1/d}}.
$$
Since $k!=\Gamma(k+1)$, this establishes (\ref{eq:vvRn}) with $C=c^{-1/d}d^{-1}$. As for the second claim of the lemma, just note that, since $\Gamma(k+\frac{1}{d})/\Gamma(k+1)$ is decreasing in $k$, we have that
$$
\sum_{k=0}^{n}\frac{\Gamma(k+\frac{1}{d})}{\Gamma(k+1)}\geq n\cdot\frac{\Gamma(n+\frac{1}{d})}{\Gamma(n+1)}=\frac{\Gamma(n+\frac{1}{d})}{\Gamma(n)},
$$
and, by Sterlings formula, $\Gamma(n+\frac{1}{d})/\Gamma(n)\sim n^{1/d}$ as $n\to \infty$.\hfill$\Box$\medskip

With Lemma \ref{le:vvRn} at hand it is not hard to prove that finiteness of $\E^*[T]$ requires finite moment of order $(d+1)/d$ for $D$.\medskip

\noindent \textbf{Proof of only-if direction of Theorem \ref{th:main}}. Write $\E^*[T]=\E^*\left[\E^*[T|D,R_D]\right]$, and, conditional on $D$ and $R_D$, define
$$
A=\left\{V^*\left(\frac{R_D}{3^{1/d}}\right)<\frac{D-1}{2}\right\},
$$
that is, $A$ is the event that at least half of the $D-1$ non-origin $\mathcal{P}^*$-points in the interior of $B(R_D)$ are located at distance at least $3^{-1/d}R_D$ from the origin. Since $\mathcal{P}^*$ has the same distribution as a Poisson process with an added point at the origin, it follows from properties of the Poisson process that the $D-1$ non-origin points are uniformly distributed in $B(R_D)$, implying that
$$
\E^*\left[V^*\left(\frac{R_D}{3^{1/d}}\right)\Big|D,R_D\right]=(D-1)\cdot\frac{\textrm{vol}[B(3^{-1/d}R_D)]}
{\textrm{vol}[B(R_D)]}=\frac{D-1}{3},
$$
where vol$[\cdot]$ denotes volume. Combining this with Markov's inequality, we can bound
$$
\PP^*(A|D,R_D)\geq 1-\frac{\E^*[V^*(3^{-1/d}R_D)]}{(D-1)/2}=\frac{1}{3}
$$
and hence
$$
\E^*[T|D,R_D]\geq \frac{\E^*[T|A,D,R_D]}{3}.
$$
On the event $A$, at least $(D-1)/2$ of the $D$ stubs at the origin must be connected to vertices whose distance to the origin is at least $3^{-1/d}R_D$, implying that
$$
T\geq \frac{D-1}{2}\cdot \frac{R_D}{3^{1/d}}.
$$
Hence, using Lemma \ref{le:vvRn}, we get that

\begin{eqnarray*}
\E^*[T] & \geq & \frac{1}{6\cdot 3^{1/d}}\cdot \E^*[(D-1)R_D]\\
& = & \frac{1}{6\cdot 3^{1/d}}\cdot \E^*[(D-1)\E^*[R_D|D]]\\
& \geq & \frac{C'}{6\cdot 3^{1/d}}\cdot \E^*[(D-1)D^{1/d}],
\end{eqnarray*}

\noindent which proves that $\E[D^{(d+1)/d}]<\infty$ is indeed necessary for $\E^*[T]<\infty$.\hfill$\Box$

\section{Repeated stable matching}

In this section we make a first attempt to formulate a pairing scheme. The algorithm is based on so called stable matchings -- see e.g.\ Gale and Shapely (1962) -- obtained by iteratively matching mutually closest points, and turns out to give finite mean for $T$ in $d\geq 2$ if $F$ has bounded support. For $d=1$ an alternative scheme is described.

For a general translation invariant homogeneous point process $\mathcal{R}$ with finite intensity, one algorithm for matching its points is as follows: First consider all pairs $\{x,y\}\subset[\mathcal{R}]$ of mutually closest points -- that is, all pairs $\{x,y\}$ such that $x$ is the closest point of $y$ and $y$ is the closest point of $x$ -- and match them to each other. Then remove these pairs, and apply the same procedure to the remaining points of $[\mathcal{R}]$. Repeat recursively. If $\mathcal{R}$ is such that almost surely $[\mathcal{R}]$ is non-equidistant  and has no descending chains, then this algorithm can be shown to yield an almost surely perfect matching, that is, a matching where each point is matched to exactly one other point. Here, the point set $[\mathcal{R}]$ is said to be non-equidistant if there are no distinct points $x,y,u,v\in[\mathcal{R}]$ with $|x-y|=|u-v|$, and a descending chains is an infinite sequence $\{x_i\}\subset[\mathcal{R}]$ such that $|x_i-x_{i-1}|$ is strictly decreasing. Furthermore, the obtained perfect matching is the unique stable matching in the sense of Gale and Shapely (1962).

Now consider our marked Poisson process, where each point $x\in[\mathcal{P}]$ has a random number $D_x\sim F$ of stubs attached to it. For each point $x$, number the stubs $1,\ldots ,D_x$ and say that stub number $i$ belongs to \emph{level} $i$ in the stub configuration. Furthermore, for $i\geq 1$, define $[\mathcal{P}]_i=\{x\in[\mathcal{P}]:D_x\geq i\}$. To connect the stubs, first match the points in $[\mathcal{P}]_1$ to each other using the stable matching algorithm described above, and pair the stubs on the first level among themselves according to this matching. Then take a uniformly chosen 2-coloring of the points in $[\mathcal{P}]_2$ such that points that are matched to each other in the matching of $[\mathcal{P}]_1$ get different colors. This gives rise to two sets of points $[\mathcal{P}]_{2,1}$ and $[\mathcal{P}]_{2,2}$ with different colors. Match the points in $[\mathcal{P}]_{2,1}$ to each other using the stable matching algorithm and do likewise with the points in $[\mathcal{P}]_{2,2}$. Connect the stubs on the second level accordingly. Then continue in this way. In general, the points in $[\mathcal{P}]_i$ are colored with a uniformly chosen $i$-coloring such that points that have been matched in the previous steps are assigned different colors, giving rise to point sets $[\mathcal{P}]_{i,1},\ldots, [\mathcal{P}]_{i,i}$ whose points are then matched among themselves using the stable matching algorithm. Since the distributions of the point sets $\{[\mathcal{P}]_{i,j}\}$ fulfill the conditions for the stable matching algorithm to yield a perfect matching, this procedure indeed provides a well-defined pairing of the stubs, and the coloring prevents multiple edges between vertices. We refer to this scheme as repeated stable matching with coloring (RSMC).

As for the total edge length per vertex with RSMC we have the following result.

\begin{prop}\label{prop:RSMC} For RSMC with degree distribution $F$, we have that
\begin{itemize}
\item[{\rm{(a)}}] $\E^*[T]=\infty$ for any degree distribution $F$ in $d=1$.

\item[{\rm{(b)}}] $\E^*[T]<\infty$ for $F$ with bounded support in $d\geq 2$.
\end{itemize}
\end{prop}

\noindent \textbf{Proof.} Part (a) is a direct consequence of Theorem 5.2(i) in Holroyd et al.\ (2008), where the expected matching distance for the stable matching of a one-dimensional Poisson process is shown to be infinite. As for (b), write $D=D_0$ and let $L_i$ denote the length of the edge created by the $i$:th stub at the origin (if $D<i$, then $L_i:=0$). Also write $p_i=\PP(D=i)$ and $p_{i+}=\PP(D\geq i)$. By a straightforward generalization of Theorem 5.2(ii) in Holroyd et al. (2008), we have that

\begin{equation}\label{eq:pp_beg}
\PP^*\left(L_i>r|D\geq i\right)\leq C\frac{i}{p_{i+}}r^{-d}
\end{equation}

\noindent for some constant $C=C(d)$. Indeed, the point sets $\{[\mathcal{P}]_{i,j}\}_{j=1,\ldots, i}$ are outcomes of translation invariant processes with intensity $p_{i+}/i$, satisfying the conditions for the existence of a well-defined stable matching, and, as pointed out by Holroyd et al., the proof of Theorem 5.2(ii) applies to any such process. It follows that $\E^*[L_i|D\geq i]\leq Ci/p_{i+}$ for some constant $C=C(d)$ and hence, with $u=\max\{i:p_i>0\}$, we have
$$
\E^*[T]=\E\Big[\sum_{i=1}^uL_i\Big]=\sum_{i=1}^up_{i+}\E^*[L_i|D\geq i]\leq C\sum_{i=1}^ui<\infty.
$$
\hfill$\Box$\medskip

\noindent \textbf{Remark.} Note that the proof applies also if the marking of the Poisson points with random degrees is not independent, but only translation invariant, that is, the probability law of the marked process is invariant under shifts of $\RR^d$. Indeed, Theorem 5.2(ii) can be applied to conclude (\ref{eq:pp_beg}) also in such a situation.\hfill$\Box$ \medskip

Clearly there is no reason to believe that RSMC is optimal to create short edges. Firstly, for a Poisson process, the stable matching is not optimal in terms of the matching distance -- see Holroyd et al. (2008) -- indicating that RSMC is not optimal within each level. Secondly, RSMC suffers from the obvious drawback that stubs on different levels are not connected to each other. Indeed, a stub on level $i$ must go at least to the nearest other vertex with degree larger than or equal to $i$ to be connected. By Lemma \ref{le:vvRn}, the expected distance to such a vertex is $Cp_{i+}^{-1/d}$, implying that
$$
\E^*[T]\geq \sum_{n=1}^{\infty}\left(\sum_{i=1}^nCp_{i+}^{-1/d}\right)p_n.
$$
This lower bound is infinite for instance for a power law distribution with sufficiently small exponent $\tau>1$, that is, for a distribution with $p_n=Cn^{-\tau}$. Indeed, then $p_{n+}=Cn^{-(\tau-1)}$, so that $p_{i+}\leq C(n/2)^{-(\tau-1)}$ for $i=n/2,\ldots,n$, and we get
$$
\E^*[T]\geq C\sum_{n=1}^\infty\frac{n}{2} \left(\frac{n}{2}\right)^{-(\tau-1)/d}n^{-\tau}=C\sum_{n=1}^\infty n^{-(\tau-1-\tau/d+1/d)},
$$
which is infinite for $\tau\leq (2d-1)/(d-1)$. In $d=2$ hence $\E^*[T]=\infty$ for $\tau\leq 3$ with RSMC, while Theorem \ref{th:main} stipulates the existence of a pairing scheme with $\E^*[T]<\infty$ for $\tau\geq 5/2$.

Before describing such a pairing scheme, we mention that, in $d=1$, RSMC can be modified to give finite mean for degree distributions with bounded support as follows: As before, let $[\mathcal{P}]_i=\{x\in[\mathcal{P}]:D_x\geq i\}$. To connect the stubs on level $i$, pick one of the two possible matchings of the points of $[\mathcal{P}]_i$ satisfying that, for any two matched points $x,y$ with $x<y$, the interval $(x,y)$ does not contain any points of $[\mathcal{P}]_i$. Call a point that is matched with a point to the right (left) right-oriented (left-oriented) and connect the stubs so that the stub at a right-oriented point $x$ is paired with the stub at the $i$:th left-oriented point located to the right of $x$. This is clearly stationary, and, since $[\mathcal{P}]_i\subset[\mathcal{P}]_{i-1}$, no multiple edges are created. We refer to this procedure as shifted adjacent matching (SAM).

\begin{prop}\label{prop:SAM}
Assume that $F$ has bounded support and let $u=\max\{p_i:p_i>0\}$. For SAM, we have that $\E^*[T]=u^2$. In particular, $\E^*[T]<\infty$ if and only if $F$ has bounded support.
\end{prop}

\noindent \textbf{Proof.} Recall that $L_i$ is the length of the edge created by the $i$:th stub at the origin ($L_i:=0$ if $D<i$). With SAM, the stub at level $i$ at the origin (if such a stub exists) will be connected to the $(2i-1)$:th point in $[\mathcal{P}]_i$ to the right or left of the origin with probability 1/2 respectively. The distance from the origin to the $(2i-1)$:th point to the right with degree at least $i$ is a sum of a NegBin$(2i-1,p_{i+})$ number of Exp(1)-variables. Hence $\E^*[L_i|D\geq i]=(2i-1)/p_{i+}$ and the proposition follows from a calculation analogous to the one in the proof of Proposition \ref{prop:RSMC}.\hfill$\Box$\medskip

\noindent \textbf{Remark.} SAM for a stub configuration on $\Z^d$ is analyzed in Deijfen and Jonasson (2006: Section 2), where it is shown that $T$ has finite mean also for \emph{stationary} degrees $\{D_i\}$ with bounded support.\hfill$\Box$

\section{Finite mean is possible when $\E[D^{(d+1)/d}]<\infty$}

We now describe a scheme that gives finite mean for the total edge length per vertex as soon as $F$ has finite mean of order $(d+1)/d$. Roughly, the model is a truncated version of RSMC, where vertices with very high degree are connected nearby (instead of having to go to other vertices with equally high degree) and the remaining stubs are then connected according to RSMC. The model is designed to make it possible to exploit results from $\Z^d$ and the procedure for connecting the high-degree vertices will involve associating the Poisson points with their nearest vertex in a uniform translation of $\Z^d$.

The following stepwise algorithm for connecting vertices with large degree to vertices with small degree in an iid stub configuration $\{\widetilde{D}_z\}$ on $\Z^d$ was formulated in Jonasson (2007, Section 3.2). It is a generalization of a similar algorithm for $\Z$ from Deijfen and Jonasson (2006) inspired by the ''stable marriage of Poisson and Lebesgue`` from Hoffman et al.\ (2006), and it is applicable also to more general graphs. For a large integer $m$, let $\widetilde{D}_z'=\widetilde{D}_z \mathbf{1}\{\widetilde{D}_z>m\}$ and say that $z$ is high if $\widetilde{D}_z'>0$, otherwise $z$ is called low. First, the positions of the vertices of $\Z^d$ are disturbed by moving each vertex independently a Unif(0,0.1)-distributed distance along a randomly chosen incident edge (this is just to make the vertex set non-equidistant). Then, in the first step, every high vertex $z$ claims its $\widetilde{D}_z$ nearest low neighbors, and a low vertex that is claimed by at least one high vertex is connected to the nearest high vertex that has claimed it. Let $\widetilde{D}_z(1)$ be the number of remaining stubs of the high vertex $z$ after this has been done. In the second step, each high vertex $z$ with $\widetilde{D}_z(1)>0$ claims its $\widetilde{D}_z(1)$ nearest low vertices that have not yet been connected to any high vertex, and each claimed low vertex is connected to the nearest high vertex that has claimed it. This is then repeated recursively. In Jonasson (2007) it is shown that, if $m$ is chosen large enough, this procedure leads to a well-defined configuration with $\E[\widetilde{T}]<\infty$ as soon as $\E[\widetilde{D}_z^{(d+1)/d}]<\infty$, where $\widetilde{T}$ denotes the total edge length (in the $\Z^d$-metric) of the origin.

Now return to the Poisson setting. For $x\in\RR^d$, write $U_x$ for the unit cube centered at $x$. Pick a point $x_0$ uniformly in the origin cube $U_{0}$ and let $\Z^d(x_0)$ be a translation of $\Z^d$ obtained by moving the origin to $x_0$. The stubs at a Poisson point $x\in[\mathcal{P}]$ are said to be associated with the point $z\in\Z^d(x_0)$ such that $x\in U_z$. Write $\widetilde{D}_z$ for the number of stubs associated with a point $z\in\Z^d(x_0)$ and $N_z$ for the number of Poisson points in $U_z$. Then $N_z\sim$ Po(1) and $\widetilde{D}_z=D_1+\ldots+D_{N_z}$, where $\{D_i\}$ are the iid marks of the Poisson points in $U_z$. The variables $\{\widetilde{D}_z\}$ are also iid and can be thought of as a stub configuration on $\Z^d(x_0)$. We number the stubs associated with each vertex $z\in\Z^d(x_0)$ randomly in some way, for instance by first numbering the Poisson points in $U_z$ randomly $1,\ldots,N_z$, and then consecutively picking one stub from each Poisson point according to that ordering until all stubs are numbered. Furthermore, in order to keep track of which Poisson point a certain stub originates from, the stubs are labelled with the position of their Poisson point. The stubs are now connected as follows:

\begin{itemize}
\item[1.] Pick $m$ large and, as in Jonasson (2007), say that a vertex $z\in\Z^d(x_0)$ is high (low) if $\widetilde{D}_z>m$ ($\widetilde{D}_z\leq m$). Match the stubs of the high vertices of $\Z^d(x_0)$ with stubs at the low vertices of $\Z^d(x_0)$ according to the algorithm from Jonasson described above, using the stubs at each vertex in numerical order. When two stubs are matched, we create an edge between the Poisson points that they originate from. Since there will be no multiple matchings between the same two vertices in $\Z^d(x_0)$, this will not give multiple edges between the Poisson points either.
\item[2.] After step 1, in $d=1$, the unconnected stubs -- that is, the stubs at the low vertices of $\Z(x_0)$ that have not been matched with stubs from the high vertices -- are connected with SAM. In $d=2$, we drop the association of the stubs with the vertices of $\Z^d(x_0)$ and take the unconnected stubs back to their Poisson points. These stubs are then connected with RSMC.
\end{itemize}

This pairing scheme is clearly translation invariant and will give finite mean for $T$ if $\E[D^{(d+1)/d}]<\infty$: If $D_i$ has finite moment of order $(d+1)/d$, then $\widetilde{D}_z$ has so as well. Hence, if the point $z\in\Z^d(x_0)$ to which the stubs at the origin vertex of $[\mathcal{P}^*]$ are associated is high, then it follows from Jonasson (2007) that $T$ has finite mean. If the point is low, then $T$ has finite mean by Proposition \ref{prop:RSMC} and \ref{prop:SAM} and the following remarks. Indeed, the remaining stubs after step 1 induce a translation invariant assignment of degrees to the Poisson points -- or, in $d=1$, to the points of $\Z(x_0)$ -- where the degrees are bounded by $m$.

\section{Further work}

We have formulated a necessary and sufficient condition for the existence of a random graph on a spatial Poisson process with finite expected total edge length per vertex. There are several related problems that remain to investigate.\medskip

\noindent \textbf{Generalization of stable matching.} A natural model for pairing the stubs is the following: Consider the marked Poisson process and, at time 0, start growing a number of balls from each Poisson point linearly in time, the number of balls that start growing from a given point being given by the mark of the point. When the collection of balls from two points meet, one ball from each point is annihilated, and an edge between the two points is created. The remaining balls continue growing (this is to avoid multiple edges). When $D\equiv 1$, so that each stub has exactly one stub attached to it, this gives rise to the unique stable matching of the points; see Holroyd and Peres (2003, Section 4). It remains to analyze the algorithm for other degree distributions, which seems to be more complicated.\medskip

\noindent \textbf{Percolation.} Apart from the edge length, other features of the configurations arising from different pairing schemes could also be investigated. One such feature is the component size, that is, the number of vertices in a given component of the graph. Will the resulting edge configuration percolate in the sense that it contains an unbounded component? This question is not relevant when $D\leq 1$ -- indeed, the configuration will then consist only of isolated edges, implying that the component size is at most 2  -- but arises naturally for other degree distributions. How does the answer depend on the degree distribution? For a given degree distribution, does the answer depend on the pairing scheme? Is it always possible to achieve percolation by taking $d$ sufficiently large?\medskip

\noindent \textbf{Independent Poisson processes.} Instead of a single Poisson process, the vertex set could be generated by two independent Poisson processes, representing two different types of vertices. This is related to matchings of points of two independent Poisson processes considered in Holroyd et al.\ (2008). We look for ways of obtaining a graph with edges running only between different types of vertices and with prescribed degree distributions for both vertex types. Can this be done if the degree distributions are different for the two vertex types? If so, how different are the degree distributions allowed to be? Which properties do the resulting configurations have?\medskip

\begin{center}
-----------------------------------------
\end{center}

\noindent \textbf{Acknowledgement} I thank Ronald Meester for discussions on the problem during a stay at VU Amsterdam spring 2008.

\section*{References}

\noindent Bollob\'{a}s, B., Riordan, O., Spencer, J. and Tusn\'{a}dy, G. (2001): The degree sequence of a scale-free random graph process, \emph{Rand. Struct. Alg.} \textbf{18}, 279-290.\medskip

\noindent Bollob\'{a}s, B., Janson, S. and Riordan, O. (2006): The phase transition in inhomogeneous random graphs, \emph{Rand. Struct. Alg.} \textbf{31}, 3-122.\medskip

\noindent Chung, F. and Lu, L. (2002:1): Connected components in random graphs with given degrees sequences, \emph{Ann. Comb.} \textbf{6}, 125-145.\medskip

\noindent Chung, F. and Lu, L. (2002:2): The average distances in random graphs with given expected degrees, \emph{Proc. Natl. Acad. Sci.} \textbf{99}, 15879-15882.\medskip

\noindent Daley, D.J. and Vere-Jones, D. (2002): \emph{An introduction to the theory of point processes}, 2nd edition, vol. I, Springer.\medskip

\noindent Deijfen, M. and Jonasson, J. (2006): Stationary random graphs on $\Z$ with prescribed iid degrees and finite mean connections, \emph{Electr. Comm. Probab.} \textbf{11}, 336-346.\medskip

\noindent Deijfen, M. and Meester, R. (2006): Generating stationary random graphs on $\mathbb{Z}$ with prescribed iid degrees, \emph{Adv. Appl. Probab.} \textbf{38}, 287-298.\medskip

\noindent Durrett, R. (2006): \emph{Random graph dynamics}, Cambridge University Press.\medskip

\noindent Gale, D. and Shapely, L. (1962): College admissions and stability of marriage, \emph{Amer. Math. Monthly} \textbf{69}, 9-15.\medskip

\noindent Hoffman, C., Holroyd, A. and Peres, Y. (2006): A stable marriage of Poisson and Lebesgue, \emph{Ann. Probab.} \textbf{34}, 1241-1272.\medskip

\noindent Holroyd, A., Pemantle, R., Peres, Y. and Schramm, O. (2008): Poisson matching, \emph{Ann. Inst. Henri Poincare}, to appear.\medskip

\noindent Holroyd, A. and Peres, Y. (2003): Trees and matchings from point processes, \emph{Elect. Comm. Probab.} \textbf{8}, 17-27.\medskip

\noindent Jonasson, J. (2007): Invariant random graphs with iid degrees in a general geology, \emph{Probab. Th. Rel. Fields}, to appear.\medskip

\noindent Kallenberg, O. (1997): \emph{Foundations of modern probability}, Springer.\medskip

\noindent Molloy, M. and Reed, B. (1995): A critical point for random graphs with a given degree sequence, \emph{Rand. Struct. Alg.} \textbf{6}, 161-179.\medskip

\noindent  Molloy, M. and Reed, B. (1998): The size of the giant component of a random graphs with a given degree sequence, \emph{Comb. Probab. Comput.}\ \textbf{7}, 295-305.\medskip

\end{document}